# DEGREE OF NEGATION
# OF AN AXIOM

Florentin Smarandache, Ph D
Professor of Mathematics
Chair of Department of Math & Sciences
University of New Mexico
200 College Road
Gallup, NM 87301, USA
E-mail: smarand@unm.edu

**Abstract.**
In this article we present the two *classical negations* of Euclid's Fifth Postulate (done by Lobachevski-Bolyai-Gauss, and respectively by Riemann), and in addition of these we propose a *partial negation* (or a degree of negation) of an axiom in geometry.

1. **Introduction.**

The most important contribution of this article is the introduction of the degree of negation (or partial negation) of an axiom and, more general, of a scientific or humanistic proposition (theorem, lemma, etc.) in any field - which works somehow like the negation in fuzzy logic (with a degree of truth, and a degree of falsehood) or like the negation in neutrosophic logic [with a degree of truth, a degree of falsehood, and a degree of neutrality (i.e. neither truth nor falsehood, but unknown, ambiguous, indeterminate)].

2. **Euclid's fifth Postulate.**

The Euclid's Fifth Postulate is formulated as follows: if a straight line, which intersects two straight lines, form interior angles on the same side, smaller than two right angles, then these straight lines, extended to infinite, will intersect on the side where the interior angles are less than two right angles.
This postulate is better known under the following formulation: through an exterior point of a straight line one can construct one and only one parallel to the given straight line.
The Euclid's V postulate (323 BC - 283 BC) is worldwide known, logically consistent in itself, but also along with other four postulates with which to form a consistent axiomatic system.
The question, which has been posted since antiquity, is if the fifth postulate is dependent of the first four?

3. **The Axiomatic System.**

An axiomatic system, in a classical vision, must be:
> 1) *Consistent* (the axioms should not contradict each other: that is some of them to affirm something, and others the opposite);
> 2) *Independent* (an axiom must not be a consequence of the others by applying certain rules, theorems, lemmas, methods valid in that system; if an axiom is proved to be dependent (results) of the others, it is eliminated from that system; the system must be minimal);
3) *Complete* (the axioms must develop the complete theory, not only parts of it).

## 4. Non-Euclidean Geometries.

The geometers thought that the V postulate (= axiom) is a consequence of the
Euclid's first four postulates. Euclid himself invited others in this research. Therefore, the system proposed by Euclid, which created the foundation of classical geometry, seemed not be independent.
In this case, the V postulate could be eliminated, without disturbing at all the geometry's development.
There were numerous tentative to "proof" this "dependency", obviously unsuccessful. Therefore, the V postulate has a historic significance because many mathematicians studied it.

Then, ideas revolved around negating the V postulate, and the construction of an axiomatic system from the first four unchanged Euclidean postulates plus the negation of the fifth postulate. It has been observed that there could be obtained different geometries which are bizarre, strange, and apparently not connected with the reality.

### 4.1. Lobachevski Geometry.

Lobachevski (1793-1856), Russian mathematician, was first to negate as
follows: "Through an exterior point to a straight line we can construct an infinite number of parallels to that straight line", and it has been named *Lobachevski geometry* or *hyperbolic geometry*. This negation is 100%.
After him, independently, the same thing was done by Bolyai (1802-1860), Hungarian from Transylvania, and Gauss (1777-1855), German. But Lobachevski was first to publish his article.
Beltrami (1835-1900), Italian, found a model (= geometric construction and conventions in defining the notions of space, straight line, parallelism) of the hyperbolic geometry, that constituted a progress and assigning an important role to it. Analogously, the French mathematician Poincaré (1854-1912).

### 4.2. Riemann Geometry.

Riemann (1826-1866), German, formulated another negation: "Through an exterior point of a straight line one cannot construct any parallel to the given straight line", which has been named *Riemann geometry* or *elliptic geometry*. This negation is also 100%.

### 4.3. Smarandache Geometries.

Smarandache (b. 1954) partially negated the V postulate (1969): "There exist straight lines and exterior points to them such that from those exterior points one can construct to the given straight lines:
   1. only one parallel – in a certain zone of the geometric space [therefore, here functions the Euclidean geometry];
   2. more parallels, but in a finite number – in another space zone;
   3. an infinite number of parallels, but numerable – in another zone of the

space;
4. an infinite number of parallels, but non-numerable – in another zone of the space [therefore, here functions Lobachevski's geometry];
5. no parallel – in another zone of the space [therefore, here functions the Riemannian geometry];
[11], [12].

Therefore, the whole space is divided in five regions (zones), and each zone functions differently. This negation is not 100% as in the *Smarandache geometries* or *mixed Non-Euclidean geometries*.
I was a student at that time; the idea came to me in 1969. Why? Because I observed that in practice the spaces are not pure, homogeneous, but a mixture of different structures. In this way I united the three (Euclidean, hyperbolic, and elliptic) geometries connected by the V postulate, and I even extended them (with other two adjacent zones).
The problem was: how to connect a point from one zone, with a point from another different zone (how crossing the "frontiers")?

In "Bulletin of Pure and Applied Science" ([3], [4], [5], [6]), then in the prestigious German magazine which reviews articles of mathematics "Zentralblatt für Mathematik" (Berlin), there exist many variants of Mixed Non-Euclidean Geometries [following the tradition: Euclid's (classical, traditional) geometry, Lobachevski's geometry, Riemann's geometry, Smarandache geometries]. Other papers were published in various journals ([1], [7], [12]), or were presented to the American Mathematical Society Meetings ([2], [8], [9]).

### 4.3.1. General Definition.
More general, a *Smarandache Geometry* is a geometry which has at least one Smarandachely denied axiom. We say that an axiom (any axiom, in any field) is *Smarandachely denied* if the axiom behaves in at least two different ways within the same space (i.e., validated and invalided, or only invalidated but in multiple distinct ways).

### 4.3.2. Example of a model of a Smarandache Geometry.
S. Bhattacharya [13] presented a simple model for a such geometry and invited the reader, as a recreational mathematics, to compose other models.

Let's consider a square ABCD and its interior points as a geometric space. A point in this space is the ordinary point, while a line is consider any segment of line that connects two opposite sides of the square. Two lines are considered parallel if they do not intersect.

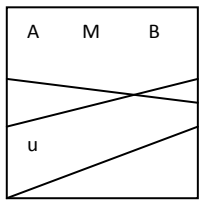

*Table 1: A Smarandache geometry on the square.*

This is a Smarandache geometry since it is partially hyperbolic non-Euclidean, partially Euclidean, and partially elliptic non-Euclidean.

Let's take a line CE and an exterior point N to it, there is an infinity of lines passing through N and parallel to CE [all lines passing through N and in between the lines (u) and (v)] – this is the hyperbolic case. But taking another exterior point M∈AB, then there is only one line parallel to CE, line AB, since only one line passes through the point M – this is the Euclidean case. Now, taking another exterior point, D, there is no parallel line passing though D and parallel to CE since all lines passing through D intersect CE – this is the elliptic case.

So, the Fifth Euclidean Postulate has been validated, but also twice invalidated.

5. **Application Quantum Mechanics.**
Ion Patrascu ([14]) proposed a model of a Smarandache Geometry, applied in quantum mechanics, built in the following way:
- an Euclidean plane α, where through any exterior point to a given line (d) there is only one parallel line;

- and an Elliptic sphere (S), where lines are defined as the big sphere circles, and points are the regular points on the sphere's surface; this is a Riemannian model of an Elliptic Geometry;

- suppose the plane \alpha cuts the sphere (S) upon a big sphere circle (C) into two equal parts; let's A and B be two distinct points on (C), which simultaneously belongs to both: the Euclidean plane α and to the Non-Euclidean sphere (S); therefore, the plane α together with the sphere (S) form a model (M) of a Smarandache Geometry. This model can be interpreted in Quantum Mechanics as follows:

- *a particle (P) that it is and it is not in a place in the same time*, is like this circle (C) which is a line [if (C) is referred to the sphere (S)] and it is not a line [if (C) is referred to the plane α] in the model (M) simultaneously;

- *a particle (R) which is in two places in the same time*, is like line AB (i.e. the line which passes through the above distinct points A and B) in the model (M); which means that 'line' AB is a straight line in the classical sense in the Euclidean plane α, while 'line' AB is the big sphere circle (C) in the Non-Euclidean sphere (S), therefore line AB is simultaneously in two different places (and has two different forms).

6. **Conclusion.**
The most important contribution of **Smarandache geometries** was the introduction of the **degree of negation of an axiom** (and more general the degree of negation of a theorem, lemma, scientific or humanistic proposition) which works somehow like the negation in fuzzy logic (with a degree of truth, and a degree of falsehood) or more general like the negation in neutrosophic logic (with a degree of truth, a degree of falsehood, and a degree of neutrality (neither true nor false, but unknown, ambiguous, indeterminate) [not only Enclid's geometrical axioms, but any scientific or humanistic proposition in any field] or partial negation of an axiom (and, in general, partial negation of a scientific or humanistic proposition in any field).
These geometries connect many geometrical spaces with different structures into a heterogeneous multi-space.
The motivation of introducing the degree of negation of an axiom is its connection with our quotidian life where the spaces are homogeneous, so an axiom or in general any proposition behaves differently in each subspace of a multi-space.

**References**

[1] Ashbacher, Charles – Smarandache Geometries – Smarandache Notions Journal, Vol. 8, No. 1-2-3, pp. 212-215, Fall 1997.
[2] Brown, Jerry L. – The Smarandache Counter-Projective Geometry – Abstracts of Papers Presented to the American Mathematical Society Meetings, Vol 17, No. 3, Issue 105, 595, 1996.
[3] Chimienti, Sandy P., Bencze, Mihály – Smarandache Anti-Geometry - Bulletin of Pure and Applied Sciences, Dehli, India, Vol. 17E, No. 1, pp. 103-114, 1998.
[4] Chimienti, Sandy P., Bencze, Mihály – Smarandache Counter Projective Geometry – Bulletin of Pure and Applied Sciences, Dehli, India, Vol. 17E, No. 1, pp. 117-118, 1998.
[5] Chimienti, Sandy P., Bencze, Mihaly – Smarandache Non-Geometry – Bulletin of Pure and Applied Sciences, Delhi, India, Vol. 17E, No. 1, pp. 115-116, 1998.
[6] Chimienti, Sandy P., Bencze, Mihály – Smarandache Paradoxist Geometry – Bulletin of Pure and Applied Sciences, Delhi, India, Vol. 17E, No. 1, pp. 123-124, 1998.
[7] Mudge, Mike – A Paradoxist Mathematician, His Function, Paradoxist Geometry, and Class of Paradoxes – Smarandache Notions Journal, Vol. 7, No. 1-2-3, August 1996, pp. 127-129; reviewed by David E. Zitarelli, Historia Mathematica, USA, Vol. 24, No. 1, p. 114, #24.1.119, 1997.
[8] Popescu, Marian – A Model for the Smarandache Paradoxist Geometry – Abstracts of Papers Presented to the American Mathematical Society Meetings, Vol. 17, No. 1, Issue 103, 265, 1996.
[9] Popov, M. R., - The Smarandache Non-Geometry – Abstracts of Papers Presented to the American Mathematical Society Meetings, Vol. 17, No. 3, Issue 105, p. 595, 1996.
[10] Smarandache, Florentin – "Paradoxist Mathematics" (1969), in Collected Papers – Vol. II, State University of Moldova Press, Kishinev, pp. 5-28, 1997.
[11] Smarandache, Florentin – Paradoxist Mathematics – Lecture, Bloomsburg University, Mathematics Department, PA, USA, November, 1995.
[12] Torretti, Roberto – A model for the Smarandache's Anti-Geometry – Universidad de Chile, Santiago, International Journal of Social Economics,

http://www.gallup.unm.edu/~smarandache/Torretti.htm
http://www.gallup.unm.edu/~smarandache/AntiGeomTorretti.pdf
[13] S. Bhattacharya, A Model to A Smarandache Geometry, J. Rec. Math., 2007.
[14] Ion Pătraşcu, *A Model of Smarandache Geometry in Quantum Mechanics*, Joint Fall 2010 Meeting of the American Physical Society Ohio Section and AAPT Appalachian and Southern Ohio Sections, Marietta College, Marietta, OH, USA, 8-9 October, 2010.